\documentclass[ 12pt,psamsfonts]{article}

\usepackage{theorem}
\usepackage{amscd}
\usepackage{babel}
\usepackage[leqno]{amsmath}
\usepackage{amssymb}
\usepackage{amsfonts}

\theoremstyle{change}
\theorembodyfont{\normalfont}
\newtheorem{atom}       {}      [subsection]

\theoremheaderfont{\normalfont\bfseries\rmfamily}
\theorembodyfont{\normalfont\itshape}

\newcommand{\ZZ}{\mathbb{Z}}
\newcommand{\QQ}{\mathbb{Q}}

\newcommand  {\shM}     {\mathcal{M}}

\newcommand  {\shL}     {\mathcal{L}}

\newcommand  {\foX}     {\mathfrak{X}}

\newcommand  {\Div}     {\operatorname{Div}}

\newcommand  {\lla}     {\longleftarrow}

\newcommand  {\lra}     {\longrightarrow}

\newcommand  {\N}       {\operatorname{N}}

\newcommand  {\NS}      {\operatorname{NS}}

\renewcommand{\O}       {\mathcal{O}}

\newcommand  {\Pic}     {\operatorname{Pic}}

\newcommand  {\pr}      {\operatorname{pr}}

\newcommand  {\quadand} {\quad\text{and}\quad}

\newcommand  {\ra}      {\rightarrow}

\newcommand  {\Sym}     {\operatorname{Sym}}

\begin{document}

\setcounter{section}{0}
\setlength{\unitlength}{1ex}

\begin{titlepage}
\title {On non-projective normal surfaces}
\author{Stefan Schr\"oer}
\end   {titlepage}


\setcounter{subsection}{0}
\maketitle

\renewcommand{\thefootnote}{}
\footnote{\hspace{-4ex}{\itshape Key words:} non-projective
surface, Neron-Severi group.}
\footnote{\hspace{-4ex}{\itshape Mathematics subject 
classification (1991)}:  14J05, 14J17.}

\begin{abstract}
In this note we construct   examples of 
{\itshape non-projective}
normal proper algebraic surfaces  and discuss the somewhat
pathological behaviour of their Neron-Severi group. 
Our surfaces   are birational to  the product of a projective
line and a curve of higher genus.  
\end{abstract}


\subsection{Introduction}

The aim of this note is to construct some simple examples of
{\itshape non-projective} normal surfaces, and discuss the
degeneration of the Neron-Severi group and its intersection form.
Here the word  {\itshape surface} refers to a 2-dimensional 
proper algebraic scheme.

The criterion of Zariski 
\cite[Cor.\ 4, p.\ 328]{Kleiman 1966} tells us that a normal
surface 
$Z $ is projective if and only if the set of points 
$z\in Z $ whose local ring 
$\O_{Z,z} $ is not 
$\QQ $-factorial allows an {\itshape affine} open neighborhood.
In particular, every resolution of singularities 
$X\ra Z $ is projective.
In order to construct 
$Z $, we therefore have to  start with a regular surface 
$X $ and  contract at least {\itshape two} suitable connected
curves 
$R_i\subset X $.

Our surfaces will be modifications of 
$Y=P^1\times C $, where 
$C $ is a smooth curve of genus 
$g>0 $; the modifications will replace some fibres 
$F_i\subset Y $ over 
$P^1 $ with {\itshape rational} curves, thereby
introducing non-rational singularities and turning  lots of
Cartier divisors into Weil divisors.

\bigskip
The Neron-Severi group 
$\NS(Z)=\Pic(Z)/\Pic^\circ(Z) $ of a non-projective surface might
become rather small, and its  intersection form might degenerate.
Our   first example  has 
$\NS(Z)=\ZZ $ and trivial intersection form. Our second  example
even has 
$\Pic(Z)=0 $. 
Our third example  allows a
birational morphism 
$Z\ra S $ to a projective   surface.

These examples provide answers for two questions concerning
surfaces posed by Kleiman
\cite[XII, rem.~7.2, p.~664]{SGA 6}. He has asked whether or not
the intersection form on the   group
$\N(X)=\Pic(X)/\Pic^\tau(X) $  of numerical classes is always
non-degenerate, and the first example shows that the answer is
negative. Here
$\Pic^\tau(X) $ is the subgroup of all invertible sheaves
$\shL$ with 
$\deg(\shL_A) = 0 $ for all curves 
$A\subset X $.  He
also has asked whether or not a normal surface with an invertible
sheaf 
$\shL $ satisfying 
$c_1^2(\shL)>0 $ is necessarily projective, and the third example
gives a negative answer.  This  should be
compared with a result on smooth complex {\itshape analytic}
surfaces
\cite[IV, 5.2, p.~126]{Barth; Peters; Van de Ven 1984}, which says
that such a surface allowing an  invertible sheaf with 
$c_1^2(\shL)>0  $ is necessarily a projective {\itshape scheme}.

\bigskip
In the following we will work over an arbitrary ground
field  
$k $  with {\itshape uncountably} many elements.  
It is not difficult to see that a normal surface over a finite
ground field is always projective. It would be interesting to
extend our constructions  to countable fields.

\subsection{A surface without ample divisors}

In this section we will construct a normal surface 
$Z $ which is not embeddable into any projective space.
The idea is to choose a suitable smooth curve 
$C $ of genus 
$g>0 $   and 
 perform certain modifications on 
$Y=P^1\times C $ called 
{\itshape mutations}, thereby destroying many Cartier divisors.

\begin{atom}
We start by choosing a smooth   curve 
$C $  such that 
$\Pic(C)\otimes \QQ $ contains uncountably many different classes
of rational points $c\in C$.
For example, let 
$C $ be an elliptic curve with at least two rational points.  We
obtain a Galois covering 
$C\ra P^1 $ of degree 2 such that the corresponding involution 
$i:C\ra C $ interchanges the two rational points. Considering its
graph we conclude that 
$i $ has at most finitely many fixed points; since there are
uncountably many rational points on the projective line, the
set 
$C(k) $ of rational points is also uncountable.

Since the group scheme of 
$n $-torsion points in the Picard scheme 
$\Pic_{C/k} $ is finite, the torsion subgroup of 
$\Pic(C) $ must be countable. Since 
$C $ is a curve of genus 
$g>0 $, any two different rational points 
$c_1,c_2\in C $  are not linearly equivalent, 
otherwise there would be a morphism 
$C\ra P^1 $ of degree 1. We conclude that 
$\Pic(C)\otimes \QQ $ contains uncountably many classes of
rational points.
\end{atom}

\begin{atom}
We will examine the product ruled surface 
$Y=P^1\times C $, and the corresponding projections 
$\pr_1:Y\ra P^1 $ and 
$\pr_2:Y\ra C $. Let 
$y\in Y $ be a rational point, 
$f:X\ra Y $ the blow-up of this point, 
$E\subset X $ the exceptional divisor, and 
$R\subset X $ the strict transform of 
$F= \pr_1^{-1}(\pr_1(y)) $. Then we can view 
$f $ as the contraction of the curve 
$E\subset X $, and I claim that there is also a contraction
  of the
curve 
$R\subset X $. Let 
$D\subset X $ be the strict transform of 
$\pr_2^{-1}(\pr_2(y)) $ and 
$\shL=\O_{X}(D) $ the corresponding invertible sheaf. Obviously,
the restriction 
$\shL\mid D $ is relatively ample with respect to the projection
$\pr_1\circ f:X\ra P^1 $; according to 
\cite{schroeer 1998} some 
$\shL^{\otimes n} $ with 
$n>0 $ is relatively base point free, hence the homogeneous spectrum
of 
$(\pr_1\circ f)_*(\Sym\shL) $ is a normal projective surface
$Z $, and
the canonical morphism 
$g:X\ra Z $ is the contraction of 
$R $, which is the only relative curve disjoint to 
$D $. We call 
$Z$ the {\itshape mutation} of 
$Y $ with respect to the center 
$y\in Y $.
\end{atom}

\begin{atom}
We observe that the existence of the contraction 
$g:X\ra Z $ is {\itshape local} over 
$P^1 $; hence we can do the same thing simultaneously for finitely
many rational points 
$y_1,\ldots, y_n $ in pairwise different closed fibres 
$F_i=\pr_1^{-1}(\pr_1(y_i)) $. If 
$f:X\ra Y $ is the blow-up of the points 
$y_i $, and 
$E_i\subset X, R_i\subset X$ are the corresponding exceptional
curves   and strict
transforms respectively, we can construct a normal proper surface 
$Z $ and a contraction 
$g:X\ra Z $ of the union
$R=R_1\cup \ldots \cup R_n $ by patching together quasi-affine
pieces over 
$P^1 $.
Since  
$Z $ is obtained by patching, there is no reason that the
resulting proper surface should be projective. We also will call 
$Z $ the {\itshape mutation} of 
$Y $ with respect to the centers 
$y_1,\ldots, y_n $.
\end{atom}

\begin{atom}
Let us determine the effect of mutations on the Picard group.
One easily sees that the maps 
$$
H^1(C,\O_C)  \lra H^1(Y,\O_{Y})  \lra H^1(X,\O_{X}) 
$$
are bijective. Let 
$\foX $ be the formal completion of 
$X $ along 
$R=\cup R_i $; since the composition 
$$
H^1(C,\O_C)  \lra  H^1(\foX,\O_{\foX}) \lra H^1(R,\O_R) 
$$
is injective, the same holds for the map on the left. Hence the
right-hand map in the exact sequence 
$$
0 \lra H^1(Z,\O_Z)         \lra H^1(X,\O_{X})  \lra  
       H^1(\foX,\O_{\foX})   
$$
in injective,  and 
$H^1(Z,\O_Z)$ must vanish. We deduce that the group scheme
$\Pic^\circ_{Z/k} $, the connected component of the Picard scheme,
is zero. Since the Neron-Severi group of 
$Y $ is torsion free, the same holds true for    
$Z $, and we conclude  
$\Pic^\tau(Z)=0 $.
\end{atom}

\begin{atom}
\label{general center}
Now let 
$F_1,F_2\subset Y $ be two different closed fibres over rational
points of
$P^1 $ and 
$y_1\in F_1 $ a rational point. The idea is to choose a second
rational point 
$y_2\in F_2 $ in a {\itshape generic}  fashion in order to
eliminate all ample divisors on the resulting mutation. Let 
$Z' $ be the mutation with respect to 
$y_1 $. By finiteness of the base number, 
$\Pic(Z') $ is a countable group, in fact isomorphic to 
$\ZZ^2 $. On the other hand,  
$\Pic(F_2) $ is uncountable, and there is a rational
point 
$y_2\in F_2 $ such that the classes of the divisors 
$ny_2 $ in 
$\Pic(F_2) $ for 
$n\neq 0 $ are not contained in the image of 
$\Pic(Z') $.
Let 
$Z $ be the mutation of 
$Y $ with respect to the centers 
$y_1,y_2 $.

I claim that there is no ample Cartier divisor on 
$Z $. Assuming the contrary, we find an ample effective divisor 
$D\subset Z $ disjoint to the two singular points 
$z_1=g(R_1) $ and 
$z_2=g(R_2) $ of the surface. Hence the strict transform 
$D'\subset Z' $ is a divisor with 
$$
D'\cap  F_2 = \left\{ y_2 \right\},
$$
contrary to the choice of 
$y_2\in F_2 $. We conclude that 
$Z $ is a non-projective normal surface.
More precisely,  there is no divisor 
$D\in\Div(Z) $ with 
$D\cdot F>0$, where  
$F\subset Z$ is a fibre over 
$P^1 $, since otherwise 
$D+nF $ would be ample for 
$n $ sufficiently large. Hence the canonical
map 
$\Pic(P^1)\ra \Pic(Z) $ is bijective,
$\Pic(Z)/\Pic^\tau(Z)=\ZZ $ holds, and the intersection form on
$N(Z)$ is zero.
\end{atom}

\subsection{A surface without invertible sheaves}

In this  section we will construct a normal surface 
$S $ with 
$\Pic(S)=0 $. We start with 
$Y=P^1\times C $, pass to a suitable mutation 
$Z $, and obtain the desired surface as a contraction of 
$Z $.

\begin{atom}
Let 
$y_1,y_2\in Y $ be two closed points in two different closed
fibres 
$F_1,F_2\subset Y $ as in 
\ref{general center} such that the mutation with respect to the
centers 
$y_1,y_2 $ is non-projective. Let 
$y_0\in Y $ be another rational point in 
$\pr_2^{-1}(\pr_2(y_1)) $, and consider the mutation 
$$
Y \stackrel{f}{\lla} X \stackrel{g}{\lra} Z
$$
with respect to the centers 
$y_0,y_1,y_2 $. We obtain a configuration of curves
on 
$X $ with the following intersection graph: 
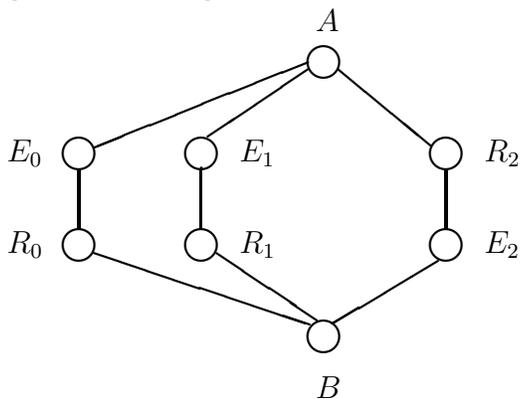
\begin{center}
\setlength{\unitlength}{0.00053300in}%
\begingroup\makeatletter\ifx\SetFigFont\undefined%
\gdef\SetFigFont#1#2#3#4#5{%
  \reset@font\fontsize{#1}{#2pt}%
  \fontfamily{#3}\fontseries{#4}\fontshape{#5}%
  \selectfont}%
\fi\endgroup%
\begin{picture}(4500,3735)(3076,-3688)
\thicklines
\put(4801,-1261){\circle{300}}
\put(3601,-1261){\circle{300}}
\put(7201,-1261){\circle{300}}
\put(3601,-2161){\circle{300}}
\put(4801,-2161){\circle{300}}
\put(7201,-2161){\circle{300}}
\put(6001,-3061){\circle{300}}
\put(5851,-436){\line(-3,-2){986.538}}
\put(5851,-361){\line(-5,-2){2094.828}}
\put(6151,-436){\line( 6,-5){900}}
\put(3601,-1411){\line( 0,-1){600}}
\put(4801,-1411){\makebox(6.6667,10.0000){\SetFigFont{10}{12}{\rmdefault}{\mddefault}{\updefault}.}}
\put(6001,-361){\circle{300}}
\put(4801,-1411){\line( 0,-1){600}}
\put(5926,-3661){\makebox(0,0)[lb]{\smash{\SetFigFont{12}{14.4}{\rmdefault}{\mddefault}{\updefault}$B$}}}
\put(7201,-1411){\line( 0,-1){600}}
\put(3751,-2236){\line( 3,-1){2115}}
\put(4951,-2236){\line( 3,-2){986.538}}
\put(7126,-2311){\line(-5,-3){1036.765}}
\put(5926,-61){\makebox(0,0)[lb]{\smash{\SetFigFont{12}{14.4}{\rmdefault}{\mddefault}{\updefault}$A$}}}
\put(2900,-1336){\makebox(0,0)[lb]{\smash{\SetFigFont{12}{14.4}{\rmdefault}{\mddefault}{\updefault}$E_0$}}}
\put(2900,-2236){\makebox(0,0)[lb]{\smash{\SetFigFont{12}{14.4}{\rmdefault}{\mddefault}{\updefault}$R_0$}}}
\put(5176,-1336){\makebox(0,0)[lb]{\smash{\SetFigFont{12}{14.4}{\rmdefault}{\mddefault}{\updefault}$E_1$}}}
\put(5176,-2236){\makebox(0,0)[lb]{\smash{\SetFigFont{12}{14.4}{\rmdefault}{\mddefault}{\updefault}$R_1$}}}
\put(7576,-1336){\makebox(0,0)[lb]{\smash{\SetFigFont{12}{14.4}{\rmdefault}{\mddefault}{\updefault}$R_2$}}}
\put(7576,-2236){\makebox(0,0)[lb]{\smash{\SetFigFont{12}{14.4}{\rmdefault}{\mddefault}{\updefault}$E_2$}}}
\end{picture}
\end{center}
Here 
$A $ is the strict transform of 
$\pr_2^{-1}(\pr_2(y_1)) $ and 
$B $ is the strict transform of 
$\pr_2^{-1}(\pr_2(y_2)) $.
Consider the effective divisor 
$D=3B+2R_0+2R_1 $; one easily calculates 
$$
D\cdot B=1, \quad D\cdot R_0=1, \quadand D\cdot R_1=1, 
$$
hence the associated invertible sheaf 
$\shL=\O_{X}(D) $ is ample on 
$D\subset X $. According to 
\cite{schroeer 1998}, the homogeneous spectrum of 
$\Gamma(X,\Sym\shL) $ yields a normal projective surface and a
contraction of 
$A\cup R_2 $. On the other hand, the curves 
$R_0 $ and 
$R_1 $ are also contractible. Since the curves 
$R_0 $, 
$R_1 $ and 
$A\cup R_2 $ are disjoint, we obtain a normal surface 
$S $ and a contraction 
$h:Z\ra S $ of 
$A $ by patching.
\end{atom}

\begin{atom}
Let 
$\shL $ be an invertible 
$\O_S $-module; then 
$\shM=h^*(\shL) $ is an invertible 
$\O_Z $-module which is trivial in a neighborhood of 
$A\subset Z $. Since the maps in
$$
\Pic(P^1)  \lra \Pic(Z)  \lra  \Pic(A)
$$
are injective, we conclude that 
$\shM $ is trivial. Hence 
$S $ is a normal surface such that
$\Pic(S)=0 $ holds.
\end{atom}

\subsection{A counterexample to a question of Kleiman}

In this section we construct a non-projective normal surface
$Z $ containing an integral Cartier divisor 
$D\subset Z $ with 
$D^2>0$. We obtain such a surface  by constructing a
non-projective normal surface 
$Z $ which allows a birational morphism 
$h:Z\ra S $ to a projective   surface 
$S $; then we can find an integral ample divisor 
$D\subset S $ disjoint to the image of the exceptional curves
$E\subset Z$.

\begin{atom}
Again we start with 
$Y=P^1\times C $ and choose two closed points 
$y_1,y_2\in Y $ as in 
\ref{general center} such that the resulting mutation is
non-projective. Let 
$y'_2 $ be the intersection of 
$F_2=\pr_1^{-1}(\pr_1(y_2)) $ with
$\pr_2^{-1}(\pr_2(y_1)) $, and 
$f:X\ra Y $ the blow up of 
$y_1 $, 
$y_2 $ and 
$y'_2 $. We obtain a configuration of curves on $X$
with intersection graph
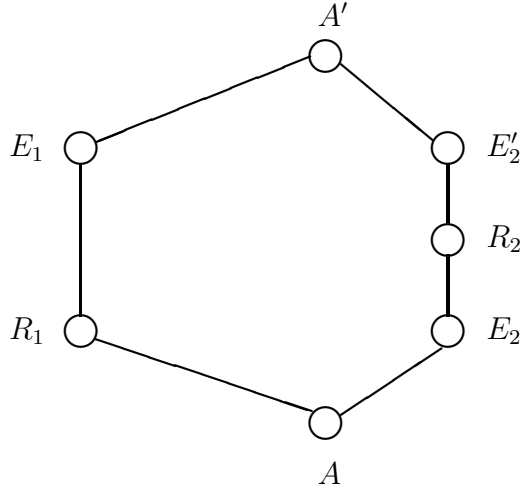
\begin{center}
\setlength{\unitlength}{0.00053300in}%
\begingroup\makeatletter\ifx\SetFigFont\undefined%
\gdef\SetFigFont#1#2#3#4#5{%
  \reset@font\fontsize{#1}{#2pt}%
  \fontfamily{#3}\fontseries{#4}\fontshape{#5}%
  \selectfont}%
\fi\endgroup%
\begin{picture}(4500,4635)(3076,-4588)
\thicklines
\put(3601,-1261){\circle{300}}
\put(7201,-1261){\circle{300}}
\put(7201,-3061){\circle{300}}
\put(6001,-3961){\circle{300}}
\put(7201,-2161){\circle{300}}
\put(3601,-3061){\circle{300}}
\put(5851,-361){\line(-5,-2){2094.828}}
\put(6151,-436){\line( 6,-5){900}}
\put(7201,-1411){\line( 0,-1){600}}
\put(6001,-361){\circle{300}}
\put(7201,-2311){\line( 0,-1){600}}
\put(7576,-3136){\makebox(0,0)[lb]{\smash{\SetFigFont{12}{14.4}{\rmdefault}{\mddefault}{\updefault}$E_2$}}}
\put(3601,-1411){\line( 0,-1){1500}}
\put(3751,-3136){\line( 3,-1){2115}}
\put(6151,-3886){\line( 3, 2){986.538}}
\put(2900,-1336){\makebox(0,0)[lb]{\smash{\SetFigFont{12}{14.4}{\rmdefault}{\mddefault}{\updefault}$E_1$}}}
\put(2900,-3136){\makebox(0,0)[lb]{\smash{\SetFigFont{12}{14.4}{\rmdefault}{\mddefault}{\updefault}$R_1$}}}
\put(5926,-61){\makebox(0,0)[lb]{\smash{\SetFigFont{12}{14.4}{\rmdefault}{\mddefault}{\updefault}$A'$}}}
\put(5926,-4561){\makebox(0,0)[lb]{\smash{\SetFigFont{12}{14.4}{\rmdefault}{\mddefault}{\updefault}$A$}}}
\put(7576,-1336){\makebox(0,0)[lb]{\smash{\SetFigFont{12}{14.4}{\rmdefault}{\mddefault}{\updefault}$E'_2$}}}
\put(7576,-2236){\makebox(0,0)[lb]{\smash{\SetFigFont{12}{14.4}{\rmdefault}{\mddefault}{\updefault}$R_2$}}}
\end{picture}
\end{center}
Here 
$A  $ is the strict transform of 
$\pr_2^{-1}(\pr_2(y_2)) $, and 
$A' $ is the strict transform of 
$\pr_2^{-1}(\pr_2(y'_2)) $.
One easily sees that there is a contraction 
$X\ra S $ of the curve
$R_1\cup R_2\cup E_2 $ and another contraction 
$X\ra Z $ of the curve
$R_1\cup R_2 $. The divisor 
$A'  $ is relatively ample on 
$S $ and shows that 
this surface is projective.
On the other hand, I claim that there is no ample divisor on 
$Z $. Assuming the contrary, we can pick an integral divisor 
$E\subset Z $ disjoint to the singularities; its strict transform 
$D\subset Y $ satisfies 
$$
D\cap F_1 = \left\{ y_1 \right\} \quadand 
D\cap F_2 = \left\{ y_2,y_2' \right\}, 
$$
where 
$F_i $ are the fibres containing 
$y_i $. Since 
$A'\cap F_2 = \left\{ y'_2 \right\} $ holds, the class of 
some multiple 
$ny_2\in \Div(F_2) $ is the restriction of an invertible 
$\O_{Y} $-module, contrary to the choice of 
$y_2 $. Hence the surface 
$Z $ and the morphism 
$Z\ra S $ are non-projective.
\end{atom}


\vspace{3em}
\noindent
        Mathematisches Institut\\
        Ruhr-Universit\"at\\
        44780 Bochum\\
        Germany\\  
        E-mail: s.schroeer@ruhr-uni-bochum.de
 
\end{document}